\documentclass[11pt]{article}
\usepackage{amssymb} \usepackage{epsf,epsfig,amsfonts}
\usepackage{amsthm} \usepackage{amsmath,amssymb,amsthm}
\usepackage{a4wide}

\newtheorem{Th}{Theorem}
\newtheorem{Lemma}{Lemma}
\newtheorem{remark}{Remark}

\newtheorem{Cor}{Corollary}
\newtheorem{example}{Example}

\newcommand{\be}{\begin{equation}}
\newcommand{\ee}{\end{equation}}

\newcommand{\rn}{\mathbb{R}^n} 

\title{Riemannian metrics having common geodesics with Berwald metrics } \date{}
\author{Vladimir S. Matveev\thanks{ Institute of Mathematics, FSU Jena, 07737 Jena Germany,  matveev@minet.uni-jena.de}}
\begin{document}
\maketitle

\begin{abstract} In Theorem 1, we generalize some  results of Szab\'o \cite{Sz1,Sz2} for Berwald metrics that  are not necessarily  strictly convex: we show that for every   Berwald metric $F$ there always exists a Riemannian metric
  affine equivalent to $F$.     As an application we show (Corollary 3) that every Berwald projectively flat metric is a Minkowski metric; this statement   is a ``Berwald" version of Hilbert's 4th  problem. 
  
  Further, we investigate  geodesic equivalence of  Berwald metrics. Theorem 2  gives a system of PDE that has  a (nontrivial) solution if and only if the given  essentially Berwald metric admits a Riemannian metric that is  (nontrivially) geodesically equivalent to it.   The system of PDE is linear and of Cauchy-Frobenius type, i.e., the derivatives of unknown functions are explicit expressions of the unknown functions.  
  As an application (Corollary 2), we obtain  that geodesic equivalence of an essentially  Berwald  metric and a Riemannian metric is always affine equivalence   provided   both metrics are complete.

\end{abstract} 
\section{Definitions and results} A 
\emph{Finsler metric} on a smooth manifold  $M$   is a function $F:TM\to \mathbb{R}_{\ge 0} $ such that:
\begin{enumerate} 
\item It is smooth on $TM\setminus TM_0$, where $TM_0$ denotes the zero section of $TM$.
\item For every $x\in M$,  the restriction $F_{|T_xM}$ is a norm on $ {T_xM}$, i.e.,  for every $\xi,\,  \eta \in T_xM$ and for every  nonnegative  $\lambda\in \mathbb{R}$  we have 
 \begin{enumerate} 
   \item $F(\lambda \cdot \xi) = \lambda\cdot   F(\xi) $, \label{2a}
   \item $F(\xi+ \eta ) \le F(\xi) + F(\eta)$, \label{2b} 
   \item $F(\xi)= 0 $ $ \Longrightarrow$ $\xi=0$.  \label{2c}
  \end{enumerate}  
  \end{enumerate} 

We always assume that   $n:= \dim(M)\ge 2$. 
We do not require  that (the restriction of) 
the function $F$ is strictly convex. In this point our definition 
is more general than the usual definition.
In addition we do not assume that the metric is reversible, i.e.,  
we do not assume that $F(-\xi)= F(\xi).$
Some  standard 
references for Finsler geometry are
\cite{Al1,BCS,BBI,Sh}.  

\begin{example}[Riemannian metric]  \label{riem}  For every Riemannian metric $g$ on $M$,  the function $F(x, \xi):= \sqrt{g_{(x)}(\xi, \xi)}$ is a Finsler metric.
\end{example}

 A Finsler metric is \emph{Berwald}, if there exists a symmetric affine connection $\Gamma$ such that the parallel transport with respect to this connection preserves the function $F$. In this case,  we call the connection $\Gamma$ the \emph{associated connection}.

  Riemannian metrics are always Berwald. For them, the associated connection coincides with the Levi-Civita connection.  We say that a Finsler metric is  \emph{ essentially Berwald}, if it is Berwald, but not Riemannian. The simplest  examples of essentially Berwald metrics are Minkowski metrics.

\begin{example}[Minkowski  metric]  \label{mink} Consider a smooth norm  on   $\mathbb{R}^n$, i.e.,  a smooth function $p:\mathbb{R}^n\to \mathbb{R}_{\ge 0}$ satisfying  \ref{2a}, \ref{2b}, \ref{2c}. We canonically identify $T\mathbb{R}^n$ with $\mathbb{R}^n\times  \mathbb{R}^n$ with coordinates $(\underbrace{x_1,...,x_n}_{x\in \mathbb{R}^n}, \underbrace{\xi^1,...,\xi^n}_{\xi\in T_x\mathbb{R}^n})$. 
Then, $F(x, \xi) := p(\xi)$ is a Finsler metric.    We see that the metric is translation invariant. Hence,  the standard flat connection  preserves  it,  i.e.,  it is a Berwald metric. If the norm  $p$ does not satisfy the parallelogram equality,  the Minkowski metric is essentially Berwald. 
\end{example}

Let $F_1,F_2$ be  Finsler metrics on the same  manifold. 
We say that $F_1$ is  \emph{geodesically equivalent}  (or  \emph{projectively equivalent})  to $F_2$, if every $F_1-$geodesic, considered as {\bf unparametrized} curve,  is also an $F_2-$geodesic. We say that they   are  \emph{affine equivalent},  if every $F_1-$geodesic, considered as {\bf parametrized} curve,  is also an  $F_2-$geodesic.  Of course, in the definition we can replace any of the Finsler metrics  by a Riemannian or  pseudo-Riemannian one, or by an affine connection. 
 
 \begin{remark} Geodesic equivalence (or affine equivalence) of Finsler metrics is not a priori a symmetric relation, as Example \ref{ex3} below shows. The reason is that for certain    Finsler metrics  the uniqueness  theorem for  the geodesics  does not hold: two different  geodesics  can have the same velocity vector, as in Example \ref{ex3} below.  Then, even  
 under the assumption that 
 all  $F_1-$geodesics are  $F_2-$geodesics, there may exist $F_2-$geodesics that are not $F_1-$geodesics. 
 
      This phenomenon evidently does not happen, if  the metrics are    strictly convex (and of course in the Riemannian case); for such metrics,  $F_1$ is  geodesically equivalent to $F_2$ if and only if $F_2$ is geodesically equivalent to $F_1$.   We will show in the beginning of Section \ref{proof} that under the assumption that the metric $F$ is Berwald, if $g$ is geodesically (or affine) 
       equivalent to $F$, then $g$ is geodesically (or  affine, resp.) equivalent to the connection $\Gamma$ associated to $F$.  
 \end{remark}
 
\begin{center} 
\begin{figure}[ht!] \label{fig}
\begin{center}{{\psfig{figure=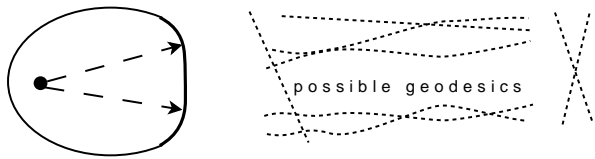,height=2cm}}}  \end{center}
  \caption{The unit sphere in the norm $p$ and possible geodesics of the corresponding Minkowski metric }
\end{figure}
\end{center} 
  \begin{example} \label{ex3} 
  Consider the Minkowski metric $F(x, \xi)= p(\xi)$ 
   such that the unit sphere $S_1:= \{\xi\in \mathbb{R}^n\mid p(\xi)=1\}$  is as on Figure \ref{fig}: the important feature of  the picture is that the part of the unit sphere lying in the marked  sector 
   is a straight line segment. Then  every curve such that its velocity vectors  are 
     in the   sector is a  geodesic. Beside such curves, the straight lines are also geodesics. We see that the standard flat metric is geodesically and affine equivalent to $F$, but the metric $F$ is neither  geodesically nor affine equivalent to the standard flat metric. 
  \end{example}

 Geodesic  equivalence  of metrics is a classical subject.  The first  non-trivial examples 
 of geodesically equivalent Riemannian metrics 
 were discovered  by   Lagrange \cite{La}. Geodesically equivalent Riemannian metrics were  studied by Beltrami \cite{Bel}, Levi-Civita \cite{LC},  Painlev\'e \cite{Pa}  and other classics.   One can find more historical details in the surveys \cite{Am,Mi} and in the introduction to the papers \cite{fomenko,Ma2}.  Geodesic equivalence of Riemannian and Finsler metrics is discussed in particular in Hilbert's 4th  problem, see \cite{Al2,Po}. Recent results on geodesic equivalence of Riemannian  and Finsler metrics  include \cite{MSB,Shen}. 

Our main results are 

\begin{Th} \label{1}  Let $F$ be a Berwald metric.
Then   there exists a Riemannian metric which is     affine equivalent to $F$. 
\end{Th} 

For strictly convex Finsler metrics, Theorem \ref{1} is due to  \cite{Sz1}. Later, other proofs were  suggested  in \cite{Sz2,To}. Our proof is similar  to the proof  in \cite{Sz2}; the modification  is based on the construction  from \cite{MRTZ}.

\begin{Th} \label{2}  Let $F$ be an essentially  Berwald metric on a connected manifold,  and  let 
$\Gamma$ be its associated connection. Suppose  a Riemannian or pseudo-Riemannian  metric  $g$ is geodesically equivalent to $F$,  but is not affine equivalent to $F$. Then there exists a constant $\mu$, a  symmetric $(2,0)-$tensor $a^{ij}$, and a  nonzero vector field  $\lambda^i$ such that the following equations are fulfilled, where ``$\ ,$'' denotes the covariant derivative with respect to  $\Gamma$:
\begin{eqnarray} 
a^{ij}_{\ \ ,k} & = &  \lambda^i\delta^j_k+  \lambda^j\delta^i_k \label{basic} \\
\lambda^i_{\ ,j} & =&  \mu \,  \delta^i_j \label{vnb}
\end{eqnarray}

  \end{Th} 

We see that  equations (\ref{basic},\ref{vnb}) are of Cauchy-Frobenius type, i.e., the derivatives of the  unknown functions   $a^{ij}$, $\lambda^i$ 
 are explicitly expressed   as functions of the unknown functions and known data (connection $\Gamma$). 
 
 \begin{remark} 
 If a Riemannian metric $g$ is affine equivalent to $F$,  equations (\ref{basic},\ref{vnb})  also have a nontrivial solution, namely   $a^{ij}= g^{ij}$, \  $\lambda^i\equiv 0$, \  $\mu =0$. 
 \end{remark} 
 
 \begin{remark} The converse of  
Theorem 2 is   also   true:  the existence of a nondegenerate 
$a^{ij}$ and of a nonzero $\lambda^i$ satisfying equations (\ref{basic},\ref{vnb}) for a certain constant $\mu$ implies the existence of a Riemannian  or a pseudo-Riemanninan metric  geodesically equivalent to $F$, but not affine equivalent to $g$. 
 \end{remark} 
 
 Recently, a system of Cauchy-Frobenius type for metrics geodesically equivalent to Berwald Finsler metrics was obtained \cite[Theorem 2]{MSB}. 
Our system  is much easier than one in \cite{MSB}: first of all, it is linear in the unknown functions, second, it contains less equations, and,  third, the equations   are much simpler  than those of \cite{MSB} and, in particular, contain no curvature terms. One cannot  obtain our equations from the equations  of \cite{MSB} by a change of unknown functions.  In order to obtain our equations from those of \cite{MSB}, one should prolong the equations of \cite{MSB} two times, and use the result of the prolongation to simplify the system.

\begin{Cor} \label{cor1}
Let $F$ be an essentially  Berwald metric on a connected closed (= compact without boundary)  manifold.  Then every Riemannian or   pseudo-Riemannian metric  geodesically equivalent to  $F$  is affine equivalent to $F$.  
\end{Cor}

\begin{Cor} \label{cor2}
Let $F$ be a complete  essentially  Berwald metric on a connected manifold. Then every  complete Riemannian or   pseudo-Riemannian metric geodesically equivalent to  $F$   is affine equivalent to $F$.
\end{Cor}

The assumptions in  Theorem \ref{2}  and 
Corollaries are important: it is possible to construct counterexamples if 
 the Berwald  metric is not essentially Berwald (i.e., is a Riemannian metric), or if one of the metrics is not complete.

\begin{Cor}[Hilbert's 4th problem for Berwald metrics] \label{cor3}  Suppose  an essentially  Berwald  metric $F$ on a connected manifold is projectively flat, that is, there exists  a   flat Riemannian 
metric  geodesically equivalent to $F$. Then $F$ is isometric to  a Minkowski metric.  
\end{Cor}

\section{Proofs}  
\subsection{Averaged  metric and proof of Theorem 1} \label{a}

Given a Finsler Berwald  metric $F$, 
we construct a Riemannian metric $g=g_F$ such that the associated connection $\Gamma$ of $F$ is the Levi-Civita connection of $g$ implying that  the metric  $g $  is affine equivalent to $F$. 
As we mentioned in the introduction,  the construction  is due to  \cite{MRTZ}, and is similar to one from  \cite{Sz2}.

Given  a smooth norm $p$  on $\mathbb{R}^{n\ge 2}$, we  canonically 
 construct a positive definite  symmetric  bilinear form $g:\mathbb{R}^n\times  \mathbb{R}^n \to \mathbb{R}$. 
 For the  Finsler metric $F$, the role of $p$ will be played by  the restriction of $F$ to $T_xM$.  
 We will see that the constructed  $g$ smoothly depends on $x$, and hence it  is a Riemannian metric.

Consider the sphere $S_1 =  \{ \xi \in \mathbb{R}^n \mid p(\xi)=1\}$.  Consider the (unique)  volume form $\Omega $ on $\mathbb{R}^n$  such that the volume of the 1-ball  $B_1 =  \{ \xi \in \rn  \mid p(\xi)\le 1\}$ is equal to $1$.
 
 Denote by $\omega$ the volume form  on $S_1$   whose value on the  vectors $\eta_1,...,\eta_{n-1}  $ tangent to $S_1$ at the point $\xi\in S_1$  is   
 given by  $\omega(\eta_1,...,\eta_{n-1}) :=   \Omega(\xi,\eta_1,\eta_2,...,\eta_{n-1})$.

Now, for every point $\xi \in S_1$, consider the symmetric  bilinear 
form $b_{(\xi)}: \mathbb{R}^n  \times \mathbb{R}^n   \to \mathbb{R} $, $b_{(\xi)} (\eta, \nu) =   D_{(\xi)}^2p^2 (\eta, \nu)$.
In this formula, $D_{(\xi)}^2p^2$ is the second differential at the point $\xi$ of the
function $p^2$ on $\mathbb{R}^n$.  The analytic expression for $b_{(\xi)} $ in the coordinates  $(\xi^1,..., \xi^n)$
  is 
  \begin{equation}  b_{(\xi)} (\eta, \nu) = \sum_{i,j} \frac{\partial^2  p^2 (\xi)}{\partial \xi^i\partial \xi^j} \eta^i \nu^j. \label{p2} 
  \end{equation}       Since the norm $p$ is  convex, the bilinear form is nonnegative definite.   Clearly, for every $\xi\in S_1 $, we have \begin{equation}   b_{(\xi)} (\xi, \xi)> 0\label{p4} \end{equation}  
    (this  is actually   the reason why we take $p^2$ and not $p$ in the definition of $b$).

Now consider  the following  symmetric bilinear $2-$form $g$ on $\mathbb{R}^n$:  for $\eta, \nu \in \rn$,  we put \begin{equation}g(\eta, \nu) = 
 \int_{S_1} b_{(\xi)}(\eta, \nu) \omega. \label{p3} 
  \end{equation} 
   We assume that the   orientation of $S_1$  is chosen in such a way that 
  $ \int_{S_1} \omega> 0$. Because of    (\ref{p4}), $g$  is positive definite.  

Now let us extend  this construction to  every tangent space $T_xM$ of the manifold, then $F_{|T_xM}$ plays
the role of $p$.  
Since the construction    depends smoothly on the point  $x\in M$, we have that $g:= g_{F}$ is a Riemannian metric on $M$.  We show that   if the metric $F$ is Berwald with the associated connection $\Gamma$, then $\Gamma$ is the Levi-Civita connection of $g$.

 Indeed, consider  a  smooth curve $\gamma$ connecting the  points $\gamma(0), \gamma(1)\in M$.  Let 
 $$\tau:T_{\gamma(0)}M\to T_{\gamma(1)}M$$ be the parallel transport of the vectors along the curve with respect to   the connection $\Gamma$. $\tau$ is a linear map. Since the metric is Berwald, $\tau$ preserves the function $F$ and, in particular, the  one-sphere $S_1$.   
 Since the forms $\Omega, \omega$ were constructed by using  the sphere $S_1$  and the linear structure of the space only,  $\tau$ 
 preserves  the form $\omega$. Since  the function $F$ is preserved   as well, everything in  formula (\ref{p3}) is preserved by the parallel transport  which implies $\tau^*g= g$. Then  $g_{ij,k}=0$, therefore   every (parametrized) geodesic of $g$ is a geodesic of $F$. Theorem \ref{1} is proved.

\subsection{Proof  of Theorem 2 and Corollaries 1, 2, 3 } 
Within  the whole section we assume that our underlying manifold is connected, orientable (otherwise we pass to an orientable cover), and has dimension at least two.  

\subsubsection{ Holonomy group of a Berwald metric $F$} 

\begin{Lemma} \label{lem1}  Let $F$ be an essentially Berwald metric on a connected manifold $M$,  and 
let $g$ be a Riemannian metric affine equivalent to $F$ (the existence of such metric is guaranteed by   Theorem \ref{1}).   
Then, the metric $g$ is symmetric of rank $\ge 2$, or   there exists one more Riemannian metric 
$h$ such that it  is  not proportional to $g$, but is affine equivalent to $g$. 
\end{Lemma} 

{\bf Proof. } We essentially 
repeat  the argumentation of \cite{Sz1,Sz2}. Take a fixed point $q\in M$. 
For  every (smooth) loop  $\gamma(t), t\in [0,1]$ with the origin in $q$ (i.e., $\gamma(0)= \gamma(1)= q$), 
we consider   the parallel transport  $\tau_{\gamma}:T_qM\to T_qM$ along  the curve. 
It is well known (see for example, \cite{Ber,Si}), that   the set 
$$
H_q:= \{  \tau_{\gamma} \mid \textrm{$\gamma:[0,1]\to M $ is a   smooth  loop, $\gamma(0)=\gamma(1)=q$}\}
$$ 
is a subgroup of the group of the orthogonal transformations of $T_qM$. Moreover, it is 
also known that  at least  one of the following conditions holds: 
\begin{enumerate} \item $H_q$ acts transitively on the unit sphere $S_1:= \{ \xi \in T_qM \mid g(\xi, \xi)=1\}$,   
\item the metric $g$ is  symmetric of rank $\ge 2$,
\item there exists one more Riemannian metric 
$h$ such that it  is  nonproportional to $g$, but is affine equivalent to $g$. 
 \end{enumerate} 
 
 In the first case, since the holonomy group preserves both $g$ and $F$, the ratio $F(\xi)^2 /g(\xi, \xi)$ is the same for all $\xi \in T_qM, \xi \ne 0$, implying that the metric $g$ is Riemannian. Lemma  \ref{lem1} is proved.

 \subsubsection{ Metrics with  degree of mobility $\ge 3$ }  \label{proof} 
 
 If the dimension of the manifold is $2$, an essentially Berwald metric is a  Minkowski metric, 
 and Theorem \ref{2} and Corollaries  \ref{cor1}, \ref{cor2}, \ref{cor3} are evident. Below, we assume that the dimension of the manifold is $\ge 3$.  
 Suppose the (Riemannian or pseudo-Riemannian) metric $\bar g$ is geodesically  equivalent to $F$, but is not affine equivalent to
 $F$.     Then  the metric $\bar g $ is  geodesically  equivalent to the averaged metric $g= g_F$, but is not affine equivalent to $g$.  If the uniqueness theorem for geodesics holds, the latter statement is trivial; for generic Finsler metrics, it probably  requires additional explanation.  
 
 In order to explain why  the metric $\bar g $ is  geodesically  equivalent to the averaged metric $g= g_F$, let us consider the set 
 $$N:= \{(x,\xi)\in TM\setminus TM_0 \ \mid \  D^2F^2_{|T_qM} \textrm{\ nondegenerate}\}.$$ 
 This  set is evidently  open. As from the following standard (see for example \cite{Ku}) argument from differential geometry it turns out, its intersection with every $T_qM \setminus TM_0  $  is not empty.

\begin{center} 
\begin{figure}[ht!] 
\begin{center}{{\psfig{figure=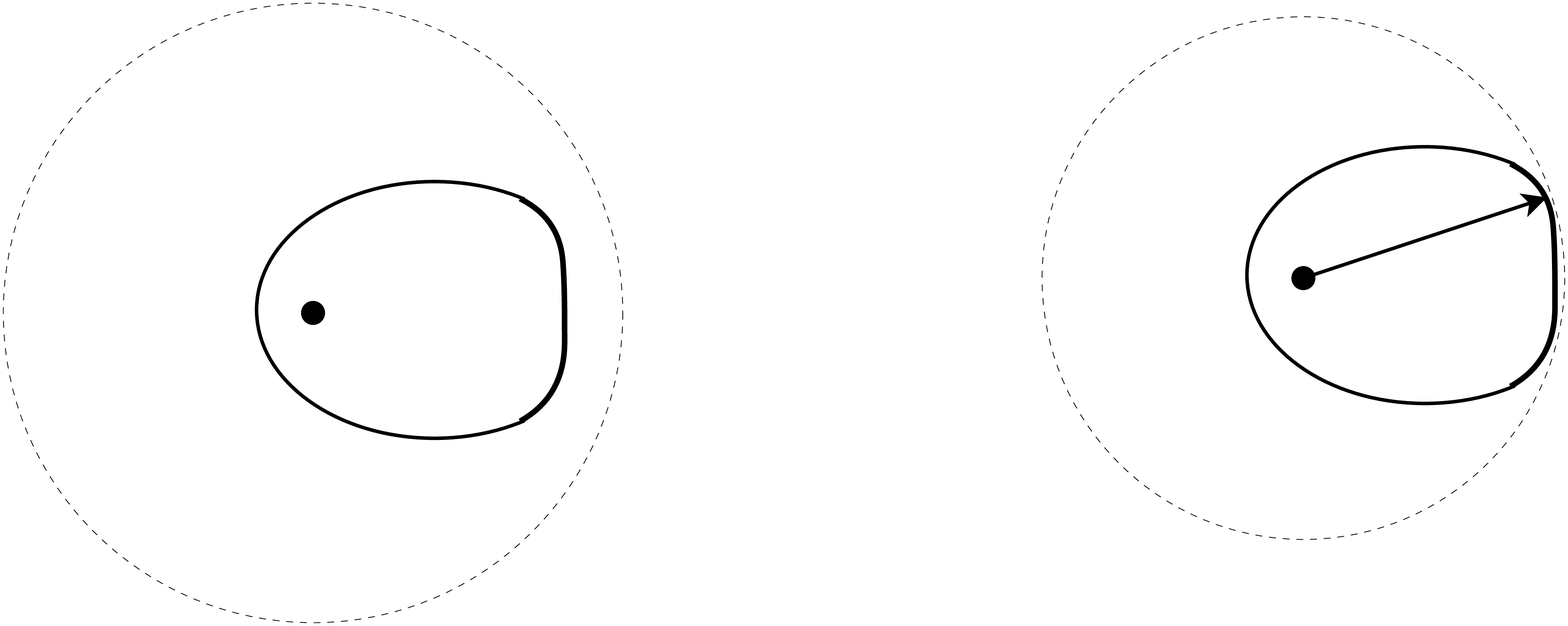,height=3cm}}}  \end{center}
  \caption{For a smooth norm $p$, there always exists a point such that the second differential of    $p^2$  is nondegenerate}\label{fig2}
\end{figure}\end{center}

 We need to show that for a smooth norm $p:= F_{|T_qM}$    on  $\mathbb{R}^n =T_qM $ there exists  a point such that  $D^2p^2$ is nondegenerate at this point. We fix an  Euclidean metric in  $\mathbb{R}^n$ and   consider the  sphere  in  $ \mathbb{R}^n$  (with respect to the chosen Euclidean metric in  $ T_qM$) 
  of large radius  such that the Finsler  sphere $S_1:= \{\xi\in T_qM \mid F(\xi)= 1\}$  lies inside, see the left-hand side of Figure \ref{fig2}.  Then, we make the radius smaller 
   until the first point of the intersection of the sphere with $S_1$, see the right-hand side of Figure \ref{fig2}. Clearly, 
   at the point of the intersection,     the second differential of    $p^2$  is nondegenerate as we claimed.

 It is well known that for $(x,\xi)\in N$ the uniqueness theorem of geodesics holds: locally, 
  there exists a   unique   $F-$geodesic $\gamma$ such that $\gamma(0)= x$ and $\dot\gamma(0)=\xi$. 
   Moreover, the geodesic $\gamma$ is also the geodesic of the associated connection $\Gamma$.    
 Then, every $\bar g-$geodesic such that $(\gamma(0),\dot\gamma(0))\in N$   is also a $\Gamma-$geodesic. Since the set $N\cap T_qM$ is open for every $q$,   the connection $\bar \Gamma$ of $\bar g$ satifies the Levi-Civita condition  
 $$\Gamma_{jk}^i-\bar \Gamma_{jk}^i - \tfrac{1}{n+1}\left(\delta_{k}^i\left(\Gamma_{j \alpha}^\alpha- \bar \Gamma^\alpha_{j\alpha}\right) + \delta_{j}^i\left(\Gamma_{k \alpha}^\alpha- \bar \Gamma^\alpha_{k\alpha}\right)\right)= 0$$  
 at every point 
 (in the proof from \cite{LC} it is sufficient to assume that only 
 the  geodesics  whose  velocity vectors are from certain open set  
 $N\subseteq TM$; $N\cap T_qM\ne \varnothing$  are common for both metrics)   implying  that  $\Gamma$ and $\bar g$  are geodesically equivalent,  and hence  $g$  and $\bar g$  are also  geodesically equivalent.

 Thus,  the metric $\bar g$   is geodesically equivalent to the averaged metric $g$ as well, 
  but not affine equivalent to $g$.  By Lemma 1, the metric $g$ is symmetric, or  there exists a Riemannian metric $h$  affine equivalent to $g$ but not proportional to $g$. We show  that if the metric $g$ is symmetric, the assumptions of Theorem \ref{1} imply that it is flat from which it follows that there exists   a metric $h= h_{ij}$  affine equivalent to $g$ but not proportional to $g$ at least on the universal cover of $M$, which is sufficient for our goals.

  By   a 
    result of 
  Sinjukov \cite{Si1}, 
   every symmetric  metric geodesically equivalent to $g$ is affine equivalent
   to $g$, unless the metric has constant curvature. 
   In the latter case, the  metric must be   flat, otherwise 
   the holonomy group discussed in the previous section  acts transitively on the unit sphere, 
   and the Finsler metric $F$  is actually Riemannian. 
   
   Thus,  at least on the universal cover of the manifold  there exists a Riemannian metric $h$  
   affine equivalent to $g$ but not proportional to $g$.

 We consider the symmetric (1,1)-tensor  
 $ a_{ij} := \left|\frac{\det(\bar g)}{\det(g)}\right|^{1/(n+1)}  
 \bar g^{\alpha \beta} g_{\alpha i} g_{\beta j}$, where $\bar g^{ij}$ is  the tensor, dual to $\bar g_{ij}$ so that $g_{i \alpha} \bar g^{\alpha j}= \delta_i^j$,  the function $\lambda:=  \frac{1}{2}a_{\alpha\beta} g^{\alpha\beta}$, and its differential $\lambda_i:= (d\lambda)_i:= \lambda_{,i}$.   
   By the result of Sinjukov \cite{Sinjukov}, 
   see also \cite{BM} and \cite{EM},  if the metric $\bar g$ is geodesically equivalent to $g$,  
 the tensor $a_{ij}$  and the (0,1) tensor  $\lambda_i$ satisfy the equation 
\begin{eqnarray} 
a_{ij ,k} & = &  \lambda_{i}g_{jk}+  \lambda_{j} g_{ik}.  \label{basic1} 
\end{eqnarray}   
  Moreover, if the metrics $g$ and $\bar g$ are not affine equivalent, $\lambda_i$ is not identically zero.

    Recall that the \emph{degree of mobility } of the metric $g$ is the dimension of the space of solutions of equation (\ref{basic1}) considered as equation on the unknown $a_{ij}$ and $\lambda_i$. In our case, the degree of mobility  is at least $3$. Indeed, $\bar a_{ij}:= g_{ij}, \bar \lambda_i:=0$ and $\hat a_{ij}:= h_{ij}, \hat \lambda_i:= 0$ are also solutions, but  by the assumptions they   are linearly independent of the solution $a_{ij}, \lambda_{i}$.

  Metrics with  degree of mobility $\ge 3$ on manifolds of dimensions $\ge 3$ were studied, in particular,  in 
  \cite{KM}, see also references therein.  The last part of the present paper will essentially use the results of  \cite{KM}, so we recommend the reader to have \cite{KM} at  hand.

 By results of \cite[Lemma 3]{KM}, 
 under the above assumptions,  for every solution $a_{ij} , \  \lambda_i$  of   equation (\ref{basic1}), in  a neighbourhood of  almost every point 
 there exists a constant $B$ and a function $\mu$  such that
 the following equations hold: 
 
\begin{eqnarray} 
\lambda_{i, j} & =&  \mu \,  g_{ij} + Ba_{ij}  \label{vnb1} \\ 
\mu_{,i} &=& 2 B \lambda_i. \label{mu}
\end{eqnarray}

 Indeed,      equation (\ref{vnb1}) 
 is  equation (30) of \cite{KM}, and  equation (\ref{mu}) is in  \cite[Remark 8]{KM} (where the function $\mu$ is denoted by $\rho$).

 Our next goal is to show that in our case $B=0$ (and,   therefore,  equations \eqref{vnb1} are fulfilled at every point of the manifold, and the  function $\mu$  is actually a  constant by   (\ref{mu})). This will also imply that     (\ref{basic1}, \ref{vnb1}) coincide with (\ref{basic},   \ref{vnb})  after raising  indices with the help of $g$.

 In order to do  this, let us  consider the solution   $A_{ij}:= a_{ij} +  h_{ij}$, $\Lambda_i :=  \lambda_i+ 0= \lambda_i$, which is the sum of the solutions $a_{ij}, \lambda_i$ and $h_{ij}, 0$.    
  The data $A_{ij}, \lambda_i$ satisfy equation \eqref{basic1}. 
  As we explained above, they therefore also satisfy  equation \eqref{vnb1} in a neighbourhood of almost every point, i.e., in a neighbourhood of almost every point  there exist a function $\tilde \mu$  and a constant $\tilde B$ 
  such that 
  \begin{equation}\label{tmp} 
  \lambda_{i,j}= \tilde \mu g_{ij} +  \tilde B (a_{ij} +  h_{ij}).  
  \end{equation} 
 Subtracting   equation \eqref{vnb1} from \eqref{tmp}, we obtain  
\begin{equation}\label{tmp1} (\mu- \tilde \mu)g_{ij} = (\tilde B- B)a_{ij} +  \tilde B h_{ij}.\end{equation}    
   We see that the right-hand side of  equation \eqref{tmp1} is a linear combination of two solution 
   $a_{ij}$ and $h_{ij}$ and is therefore also a solution of \eqref{basic1} (with an appropriate   $\lambda_i$). As it was proved in \cite[Lemma 1]{BKM} (the result is essentially due to Weyl \cite{We}), the function $\mu- \tilde \mu$ must be a constant. Since $g, a, $ and $h$ are linearly independent, all coefficients in the linear combination  \eqref{tmp1} are zero implying $B=0$.

    Thus,  equations (\ref{basic1}, \ref{vnb1}) coincide with  equations (\ref{basic},\ref{vnb}) after raising the indexes. Theorem \ref{2} is proved.

{ \bf Proof of Corollaries \ref{cor1},\ref{cor2}.}  As we explained  above, we can assume that the dimension of the manifold is $\ge 3$ and the degree of mobility is $\ge 3$. Under these assumptions, Corollary \ref{cor1} follows from      \cite[Theorem 2]{KM} 
(if $g$ is Riemannian,  the result is due to  \cite[Theorem 16]{Ma2}; in view of Theorem \ref{2}, the result follows from    \cite[Theorem 5]{Mi2}), and   Corollary \ref{cor1} follows  from  \cite[Theorem 2]{Ma1} (if $g$ is  Riemannian,  the result is due to  \cite[Theorem 1]{KM}). 
 
  {\bf Proof of Corollary \ref{cor3}.}  Suppose that a flat Riemannian metric $\bar g$ is geodesically equivalent to an
   essentially  Berwald metric $F$.   Consider the averaged metric $g= g_F$ constructed in Section  \ref{a}. It is affine equivalent to $F$, and,  therefore, as we explained in Section \ref{proof},   is geodesically  equivalent to $\bar g$.

     By the classical Beltrami Theorem  (see for example \cite{beltrami_short}, or the original papers \cite{Bel} and \cite{schur}),  the metric $g$ has constant curvature. If the curvature of $g$ is not zero, the holonony group of $g$ acts transitively on the unit sphere implying the metric $F$ is actually Riemannian.  Thus, the metric $g$ is flat. Then, there exists a coordinate system such that $\Gamma\equiv 0$. In this coordinate system, parallel transport along a curve 
     does not depend on the curve and is the usual parallel transation $x \mapsto x+ T$. Since the parallel transport preserves $F$,  we have that $F$ is translation-invariant implying it is Minkowski metric as we claimed.

  {\bf Acknowledgements. }
I   thank    Deutsche Forschungsgemeinschaft (Priority Program 1154 --- Global Differential Geometry)
  and FSU Jena for partial financial support, and V. Kiosak for attracting my attention to the paper \cite{MSB} and for useful discussions. The idea to use the construction from \cite{MRTZ} in the proof of Theorem \ref{1} appeared during the lunch  with L. Kozma and  T. Q.  Binh; I thank them and the University of Debrecen for their  hospitality.   I also thank J. Mikes and Z. Shen  for their comments, and the referee for grammatical and stylistic corrections.



\begin{thebibliography}{999}

\bibitem[Al1]{Al2} J. C. {\'A}lvarez Paiva: {\em
Symplectic geometry and Hilbert's fourth problem. } 
J. Differential Geom. {\bf 69}(2005), no. 2, 353--378. 

\bibitem[Al2]{Al1} J. C. {\'A}lvarez Paiva: {\em Some problems on Finsler geometry.} Handbook of differential geometry. Vol. II,  1--33, Elsevier/North-Holland, Amsterdam, 2006.

\bibitem[Am]{Am} A. V.  Aminova: {\it Pseudo-Riemannian manifolds with general geodesics.}   Russian Math. Surveys  48  (1993),  no. 2, 105--160, MR1239862, Zbl 0933.53002. 



\bibitem[BCS]{BCS} D. Bao, S.- S. Chern, Z. Shen: {\em An Introduction to
Riemann-Finsler Geometry.} Grad. Texts Mathem. {\bf 200} Springer Verlag
New York 2000


 \bibitem[Bel]{Bel} E. Beltrami:
{\it
Resoluzione del problema: riportari
i punti di una
superficie sopra un piano in modo che le linee geodetische vengano
rappresentante da linee rette.}
Ann. Mat., {\bf 1}(1865), no. 7, 185--204.
 
 
\bibitem[Ber]{Ber} M. Berger: {\em Sur les groupes d'holonomie homog\`ene des vari\'et\'es \`a  connexion affine et des vari\'et\'es riemanniennes.}  
Bull. Soc. Math. France {\bf 83} (1955), 279--330. 

\bibitem[BM]{BM} A.  V. Bolsinov, V. S. Matveev:
{\it Geometrical interpretation of Benenti's systems},
J. of Geometry and Physics,  {\bf 44}(2003),  489--506, MR1943174, Zbl 1010.37035.

\bibitem[BKM]{BKM}  A. V. Bolsinov, V. Kiosak, V. S. Matveev:  
    \emph{ Fubini Theorem for pseudo-Riemannian metrics.} arXiv:0806.2632. 

\bibitem[BBI]{BBI}%
{ D. Burago,  Yu. Burago, S. Ivanov: }%
{\em A course in metric geometry.}
  Graduate Studies in Mathematics, {\bf 33},
AMS, Providence, RI, 2001.

\bibitem[EM]{EM} 
M. Eastwood,  V. S. Matveev:  \emph{ Metric connections in projective differential geometry.}
 Symmetries and Overdetermined Systems of Partial Differential Equations (Minneapolis, MN, 2006), 339--351,
 IMA Vol. Math. Appl.,    {\bf
144}(2007),   Springer, New York. 	arXiv:0806.3998.


\bibitem[KM]{KM} V. Kiosak, 
V. S. Matveev:
{\em  Proof of projective Lichnerowicz conjecture for pseudo-Riemannian metrics with degree of mobility greater than two.}   arXiv:0810.0994


\bibitem[La]{La} J.-L.  Lagrange: \emph{ Sur la construction des cartes g\'eographiques.} Nouveaux M\'emoires de l'Acad\'emie des Sciences et Bell-Lettres de Berlin, 1779. 


 \bibitem[LC]{LC}
 T. Levi-Civita: {\it Sulle trasformazioni delle equazioni
 dinamiche}.  Ann. di Mat., serie $2^a$, {\bf 24}(1896), 255--300.


 \bibitem[Ku]{Ku} W. K\"uhnel: {\em 
Differential geometry. 
Curves---surfaces---manifolds. }   Student Mathematical Library, 16. American Mathematical Society, Providence, RI, 2002.


\bibitem[Ma1]{fomenko}  V. S. Matveev:  { \it  Beltrami problem, Lichnerowicz-Obata conjecture
and  applications of integrable systems in differential geometry,} Tr. Semin.
  Vektorn.   Tenzorn. Anal, {\bf 26}(2005), 214--238. 

\bibitem[Ma2]{beltrami_short}
V. S. Matveev:
\emph{Geometric explanation of Beltrami theorem}.
Int. J. Geom. Methods Mod. Phys. \textbf{3} (2006), no. 3, 623--629.

\bibitem[Ma3]{Ma1}
V. S. Matveev:
{\em  On degree of mobility of complete metrics.}  
Adv. Stud. Pure Math., {\bf 43}(2006), 221--250. 

\bibitem[Ma4]{Ma2}
V. S. Matveev:
{\em Proof of projective Lichnerowicz-Obata conjecture}.
 J. Diff. Geom.  {\bf 75}(2007),  459--502, 	arXiv:math/0407337




\bibitem[MRTZ]{MRTZ}%
{ V. S. Matveev,  H.-B. Rademacher, M. Troyanov, and A. Zeghib: }%
{\em Finsler conformal Lichnerowicz-Obata conjecture.}  Accepted to Ann. Fourier,  	arXiv:0802.3309v2

\bibitem[Mi1]{Mi2}
 J. Mikes: {\it Global geodesic mappings and their generalizations for compact Riemannian spaces. Differential geometry and its applications. }  Proceedings of the 5th international conference, Opava, Czechoslovakia, August 24-28, 1992.  Silesian Univ. Math. Publ. {\bf  1}(1993), 143-149.  http://www.emis.de/proceedings/5ICDGA/III/mike.ps

\bibitem[Mi2]{Mi}
 J. Mikes: {\it
  Geodesic mappings of affine-connected and Riemannian spaces.
Geometry, 2.}
  J. Math. Sci. {\bf 78}(1996), no. 3, 311--333. 
  
  
  
\bibitem[MBB]{MSB}%
 J. Mikes,  S. B\'acs\'o,  V. Berezovski: {\em  Geodesic mappings of weakly Berwald spaces and Berwald spaces onto Riemannian spaces.}  Int. J. Pure Appl. Math.  {\bf 45}(2008),  no. 3, 413--418. 


\bibitem[Pa]{Pa} P. Painlev\'e: \emph{Sur les int\'egrale quadratiques des
\'equations de la Dynamique}. Compt.Rend., {\bf 124}(1897),
221--224.


\bibitem[Po]{Po} A.V. Pogorelov: {\em Hilbert's Fourth Problem.}  Scripta Series in Mathematics, Winston and Sons, 1979.


\bibitem[Sh1]{Sh} Z. Shen: {\em Lectures on Finsler geometry.} World Scientific,
Singapore 2001


\bibitem[Sh2]{Shen} Z. Shen: {\em On projectively related Einstein metrics in Riemann-Finsler geometry.}   Math. Ann.  {\bf 320}(2001),  no. 4, 625--647.


\bibitem[Sc]{schur} F. Schur: \emph{ Ueber den Zusammenhang der R\"aume constanter
Riemann'schen Kr\"ummumgsmaasses mit den projektiven R\"aumen. } Math. Ann. {\bf 27}(1886), 537--567. 
\bibitem[Sim]{Si} J. Simons: {\em On transitivity of holonomy systems.}  Annals of Math.
{\bf 76}(1962) 213--234.



\bibitem[Si1]{Si1} N. S.  Sinjukov: \emph{ On geodesic mappings of Riemannian spaces onto symmetric Riemannian spaces. } (in Russian) Dokl. Akad. Nauk SSSR (N.S.) {\bf  98}(1954), 21--23.


\bibitem[Si2]{Sinjukov} N. S.  Sinjukov:
 {\it Geodesic mappings of Riemannian spaces},  (in Russian)
``Nauka'', Moscow, 1979, MR0552022, Zbl 0637.53020.

\bibitem[Sz1]{Sz1} Z. I. Szab\'o:  {\em  Positive definite Berwald spaces (Structure theorems)}.
Tensor N. S. {\bf 35}(1981) 25–39.

\bibitem[Sz2]{Sz2} Z. I. Szab\'o:
{\em Berwald metrics constructed by Chevalley's polynomials.} {\tt arXiv:math/0601522}


\bibitem[To]{To} R. G. Torrome:
{\em Average Riemannian structures associated with a Finsler structure.}
{\tt arXiv:math/0501058v5}

\bibitem[We]{We}  H. Weyl: {\em Zur Infinitisimalgeometrie: Einordnung der projektiven
und der  konformen Auffassung.} Nachrichten von der K. Gesellschaft
der Wissenschaften zu G\"ottingen, Mathematisch-Physikalische
Klasse, 1921;
 ``Selecta Hermann Weyl'', Birkh\"auser Verlag,
   Basel und Stuttgart,
1956.

\end{thebibliography}
\end{document}